\documentclass{amsart}
\usepackage{amscd}
\usepackage{latexsym}
\usepackage{amssymb,amsmath}
\usepackage{amsfonts}
\usepackage{url}

\newcommand{\Q}{{\mathbb Q}}

\renewcommand{\P}{\mathbb{P}}
\newcommand{\Proj}{{\mathop{\bf Proj}}}

\theoremstyle{plain}
\newtheorem{thm}{Theorem}[section]

\theoremstyle{definition}

\theoremstyle{remark}

\def\O{{\mathcal O}}

\begin{document}
\bibliographystyle{plain}
\bibstyle{plain}

\title{Appendix to a paper of B. Williams}

\author{Adam Logan}
\address{The Tutte Institute for Mathematics and Computation,
  P.O. Box 9703, Terminal, Ottawa, ON K1G 3Z4, Canada}
\address{School of Mathematics and Statistics, 4302 Herzberg Laboratories,
  1125 Colonel By Drive, Carleton University, Ottawa, ON K1S 5B6, Canada}
\email{adam.m.logan@gmail.com}

\date{\today}

\begin{abstract}
  In this appendix to a paper by B. Williams,
  we give birational equivalences between
  the models of the Hilbert modular surfaces for $\Q(\sqrt{29})$ and
  $\Q(\sqrt{37})$ given there and those previously found by Elkies and Kumar.
\end{abstract}

\maketitle

\section{Introduction}\label{sec:intro}
The main results of the paper \cite{williams} give explicit presentations of
the ring of Hilbert modular forms for the full ring of integers of
$\Q(\sqrt{29})$ and $\Q(\sqrt{37})$.  Applying the $\Proj$ functor to these
rings gives the {\em Baily-Borel compactification} \cite[Theorem 2.7.1]{vdg}
of the associated Hilbert
modular varieties.  However, these schemes are somewhat inconvenient to work
with in a computer algebra system, because they are embedded in a highly
singular weighted projective space of large dimension.  In addition, the
cusp singularities, though they admit elegant and meaningful resolutions
described in \cite[Section 2.6]{vdg}, are not canonical (\cite [(1.1), (1.2)]{reid}),
and therefore interfere with
calculations of linear systems on the surfaces.  

On the other hand, the paper \cite{elkies-kumar} gives many models of
Hilbert modular surfaces as subvarieties of ordinary or weighted projective
space with only canonical singularities, and as elliptic surfaces over a
rational curve.  These are very convenient for such purposes as 
calculating rational curves and elliptic fibrations on the surface, 
but they do not represent easily described functors of families of
abelian surfaces with real multiplication.  Thus it may not be clear, for
example, how to describe the Hirzebruch-Zagier divisors \cite[Chapter 5]{vdg}
on them.

To bridge the gap, therefore, it is convenient to exhibit an explicit birational
equivalence between the Baily-Borel compactifications of \cite{williams} and
the elliptic surfaces of \cite{elkies-kumar}.  We find these in two stages:
\begin{enumerate}
\item By repeatedly projecting away from bad components of the
  singular locus, we obtain a model of the surface in projective space with
  only canonical singularities.  (Recall from, for example, \cite[(1.2)]{reid}
  that canonical singularities of a surface are exactly those with a
  resolution by rational curves of self-intersection $-2$ whose intersection
  graph is a Dynkin diagram of type $A, D, E$.)
\item Having found a suitable model, we find rational curves and elliptic
  fibrations on it until we come across a fibration
  that matches the description given
  in \cite{elkies-kumar}.  This is an iterative process: given some rational
  curves, we can find the $ADE$ configurations supported on them and write
  down the corresponding elliptic fibrations; in turn, an elliptic
  fibration on a K3 surface of large Picard rank is likely to have reducible
  fibres, which are composed of rational curves.
\end{enumerate}

It is not clear how to guarantee success for either of these stages,
but in both of the examples discussed here this process could be carried
through in Magma \cite{magma}.  The code supporting the claims of this
appendix is available as \cite{magma-scripts}.

\section{Finding a model with canonical singularities}
We recall \cite[Theorem 7.3.3]{vdg} that the Hilbert modular surfaces for
$\Q(\sqrt{29})$ and $\Q(\sqrt{37})$ are birationally equivalent to K3 surfaces.
If $D$ is a nef and
big divisor on a K3 surface, then by Riemann-Roch we have $D^2 = 2h^0(D)-4$.
So if $S \subset \P^n$ is birational to a K3 surface, we can use $\deg S - 2n$
as a measure of how far we are from finding a canonical model given by a
complete linear series.

Our first step for both surfaces will be to eliminate some of the variables
to obtain a surface in the weighted projective space $\P(2,2,3,3,6)$.  This
space is embedded in $\P^7$ by $\O(6)$, so we obtain a birationally equivalent
surface in ordinary projective space.  The general way to improve a singularity
is to blow up along a component of the singular locus, which in computer
algebra can often be interpreted as a map
$\P^n \dashrightarrow \P^n \times \P(a_1, \dots, a_k)$, where the $a_i$ are the
degrees of the equations defining the component.  In a situation such as ours,
however, this would rapidly lead to an unmanageable profusion of variables.
Instead we use two basic steps: projection away from (the linear span of) a
component of the singular locus and the $2$-uple Veronese embedding.  The
disadvantage of this is that there is no guarantee that a step will actually
improve the singularities; projecting away from a point, for example, is locally
the same as blowing it up, if there are no lines through the point.  Such
lines will be contracted to new singularities.  In addition, the projection
may be a map of degree $2$ rather than a birational equivalence to its image.
We will
start by projecting away from $1$-dimensional components of the singular
subscheme and then proceed to components of large degree; although the degree
of the component of the singular locus of a scheme at a singularity of type
$A_i, D_i, E_i$ is $i$ and may be arbitrarily large in the first two cases,
each such singularity
contributes $i$ to the Picard number, which is at most $20$ for a K3 surface in
characteristic $0$.  Thus we expect that a component of the singular
locus of degree much greater than $10$ is not a canonical singularity (and
can verify this expectation using Magma if necessary).

We now give the details of the construction, first for $\Q(\sqrt{29})$ and then
for $\Q(\sqrt{37})$.  For modular forms, we use the notation of \cite{williams};
for coordinates on ordinary projective space $\P^n$, we use $x_0, \dots, x_n$.
It is never necessary to consider more than one projective space at a time, so
no confusion will result from this.

\subsection{$\Q(\sqrt{29})$}
As noted above, we begin by defining a rational map from the Baily-Borel
compactification
in the weighted projective space $\P(2,2,3,3,4,5,6,6,7,8,9)$ to a surface in
$\P(2,2,3,3,6)$ by the equations $(E_2:\phi_2:\phi_3:\psi_3:\phi_6)$.
Every other generator of the ring of modular forms can be eliminated using
a relation of degree $1$ in it (these relations can be chosen in order so as
not to reintroduce previously eliminated generators).  Making these
substitutions introduces two extraneous components to the scheme, supported
at $\phi_2 = \phi_3 = 0$ and $E_2 = \phi_2 = \phi_3 = 0$, but these
are easily removed.  Since one of the equations for the desired component
is $E_2\phi_2^2 + 4\phi_2^3 - \phi_3^2 + 4\phi_3\psi_3$ of degree $6$, the
embedding by $\O(6)$ goes naturally into $\P^6$; the image has degree $15$.

The singular subscheme of the image in $\P^6$ includes a highly nonreduced
component of degree $45$ supported at $(1:0:0:0:0:0:0)$.  We 
project away from this point to obtain a surface of degree $11$ in $\P^5$.
In turn, the singular subscheme of this surface has a component of degree
$21$ supported at $(0:0:0:0:0:1)$; projecting away from this point gives a
surface of degree $8$ in $\P^4$, which is a complete intersection.
It does not have any isolated rational singular points to project from, so
it is time to use the Veronese embedding.  We consider the map given by the
linear system of quadrics that vanish on the two singular lines
$x_0 = x_1 = x_2 = 0$ and $x_0 = x_2 = x_3 = 0$ as well as on the isolated
singular points $(\pm \frac{\sqrt{-1}}{4}:0:0:0:1)$, modulo the quadric
vanishing on the surface.  This gives a map to a surface in $\P^6$, whose
singularities are all isolated points.

The surface in $\P^6$ has a bad singularity at $(0:0:0:0:0:0:1)$, but
projection from this point is not a birational equivalence from the surface
to its image.  Instead, we
project away from the mild singularities at $(-1:0:1:0:0:0:0)$ and
$(0:0:0:1:0:0:0)$, obtaining a second surface in $\P^4$.  On this surface
we consider the linear system of quadrics vanishing on the reduced subscheme
of the singular locus, which consists of two lines defined over $\Q$
and two triples of points
defined over cubic fields (in fact the same cubic fields that will appear later
as fields of definition of the bad fibres of the elliptic fibration).
After taking the quotient by the quadric in the ideal of the surface in
$\P^4$, this gives us a map to $\P^5$.  The only bad singularity of the image
is at $(0:0:0:0:0:1)$, and projecting away from it gives a surface
$S_{29} \subset \P^4$ with four isolated double points and no other
singularities.
This surface is defined by the equations
\begin{equation}
  \label{k3-29}\begin{split}
&-48x_1^2 + 68x_1x_2 + 3x_2^2 + 20x_0x_3 + 232x_1x_3 - 42x_2x_3 - 
    181x_3^2 - 3x_4^2 = \cr
&320x_0x_1^2 + 2496x_1^3 + 880x_0x_1x_2 - 6560x_1^2x_2 - 20x_0x_2^2
    + 3628x_1x_2^2 + 114x_2^3 \cr
&    \quad - 17584x_1^2x_3 + 20124x_1x_2x_3 - 
    2441x_2^2x_3 + 29244x_1x_3^2 - 13032x_2x_3^2 \cr
&    \quad - 14117x_3^3 +  20x_0x_4^2 + 156x_1x_4^2 - 114x_2x_4^2 - 195x_3x_4^2 = 0.
\end{split}\end{equation}
Being a complete intersection of degrees $2,3$ in $\P^4$ with only canonical
singularities, it is a birational model of a K3 surface.  It is easily checked
that every step along the way was invertible, so it is birationally equivalent
to the original surface.

\subsection{$\Q(\sqrt{37})$}
This calculation is quite similar to the one just presented.  We use the
notation of the supplementary material to \cite{williams} for modular forms:
in terms of the notation in the paper, we have
\begin{equation}
  g_2 = \phi_1 \psi_1, \quad g_{3,1} = \phi_1^3, \quad g_{3,2} = \phi_1\psi_2, \quad g_{6,2} = \psi_2 \phi_4.
\end{equation}
We begin by
expressing the other generators in terms of $E_2$ and these four,
which gives a surface in $\P(2,2,3,3,6)$, and embedding that surface into
$\P^6$ (again there is an equation of degree $6$).  This surface has degree
$24$ rather than $15$ as in the previous
example, so we expect the reduction process to take a bit longer.

Now we project away from the worst isolated singularity at
$(0:0:0:0:0:0:1)$, and then from the image of the second worst at
$(1:0:0:0:0:0)$.  After that the worst isolated singularity is at
$(0:1:0:0:0)$; projection gives a surface of degree $7$ in $\P^3$.  We
have decreased the invariant $\deg S - 2n$ by $11$, but there is still
work to do.

We map back to $\P^6$ by the forms of degree $2$ vanishing at the three
worst singular points $(0:0:0:1),(0:0:1:0),(1:0:0:0)$, obtaining a surface in
$\P^6$.  This surface is singular along $x_5 = x_6 = 0$, so we project away
from that line.  From the resulting surface in $\P^4$, we map by forms of
degree $2$ vanishing on the reduced subscheme of the singular locus mod
those vanishing on the surface to get a surface of degree $6$ in $\P^3$.
In turn, we map by quadrics vanishing along the singular line $x_2 = x_3 = 0$
and the point $(-1:1:0:0)$ to a surface in $\P^5$, which has a bad singularity
at $(-1:1:-1:1:0:0)$.  Projecting away from this, we return to $\P^4$,
and from there we map by quadrics vanishing on the $1$-dimensional components
of the singular locus (whose degrees are $2, 1, 1$) mod those vanishing on the
surface to obtain a surface $S_{37}$ in $\P^5$.  This surface is defined by the
equations
\begin{equation}\label{k3-37}\begin{split}
&-x_3x_4 + x_0x_5 + 2x_1x_5 + x_2x_5 = \cr
&-243x_1^2 + 243x_0x_2 + 162x_1x_3 - 324x_2x_3 + 3x_4^2 - 8x_4x_5 +
  4x_5^2 = \cr
&2187x_0^2 + 8748x_0x_1 - 5184x_0x_3 - 6480x_1x_3 - 1296x_2x_3 \cr  
\quad& + 2592x_3^2 + 40x_4x_5 - 76x_5^2 = 0.\cr
\end{split}\end{equation}

This is a smooth complete intersection of quadrics in $\P^5$, so it is a K3
surface.  Again all maps are easily checked to be invertible.

\section{Finding rational curves and elliptic fibrations}
We now indicate how to find elliptic fibrations on the surfaces defined by
(\ref{k3-29}), (\ref{k3-37}) that match those of \cite[Sections 16, 18]{elkies-kumar}.
We do not use a systematic procedure for this; we simply list rational
curves, group them into fibres of elliptic fibrations, and hope for the best.
The basic fact (going back to Pyatetski\u\i-Shapiro and Shafarevich if not
beyond) that allows us to construct elliptic fibrations is the following:

\begin{thm} 
  Let $R_1, \dots, R_k$ be smooth rational curves on a K3 surface in the
  configuration of a reducible fibre of a minimal elliptic fibration.
  Then there is a genus $1$ fibration on the surface one of whose
  fibres is supported on the $R_i$. 
\end{thm}

Indeed, a suitable linear combination of the $R_i$ is effective and nonzero
and has self-intersection $0$ and no base components.  The result then follows
from Riemann-Roch.  (See, for example, \cite[Table 4.1]{silverman} for the
list of configurations.)

Once we know that an elliptic fibration exists, we still have to construct it.
This is quite easy on $S_{37}$ (\ref{k3-37}), because it is smooth: given a
scheme-theoretic fibre $F$, consider a polynomial $p$ of degree $d$ (to be
chosen as small as possible) that vanishes on $F$ but not on the whole surface.
Let $R = (p=0) \setminus F$ be the residual divisor; then for any form $q$
of degree $d$ vanishing on $R$, the residual $(q = 0) \setminus R$ is linearly
equivalent to $F$.  So the linear system of such polynomials modulo those
vanishing on the surface defines the fibration (we are implicitly using the
projective normality of surfaces defined by a complete linear system).  On the
singular surface $S_{29}$ (\ref{k3-29}), things could be slightly more
difficult if we had to consider fibres with multiple components that meet the
singular locus or that are supported there, but we do not so the same method
works.

After equations defining a fibration are found, it is a simple matter to write
down the fibre over the generic point of $\P^1$ and (in the cases described here
where the fibre is of low degree) to find an explicit map to a Weierstrass
model.  Using Magma we can then find the reducible fibres and check whether
they match those of \cite{elkies-kumar}.  If not, we instead decompose them
into new rational curves, form them into new fibres, and repeat the process.

\subsection{$\Q(\sqrt{29})$}
We are working on a surface with four singular points of which one is rational.
The tangent cone at the singular point meets the surface in a conic $C_1$,
two lines $L_1,L_2$, where $L_1$ contains $(0:0:-1:0:1)$ and $L_2$ contains
$(0:0:1:0:1)$,
and two conjugate lines defined over $\Q(\sqrt{-1})$.  When we project away
from the singular point, we obtain a quartic in $\P^3$ with two rational
singular points (the images of the rational lines).  Pulling back the curves
in the tangent cones at these singular points, we obtain two additional
lines $L_3,L_4$ not containing the singular point of $S_{29}$: for definiteness,
let $L_3$ be the one that contains $(1:0:-1:0:1)$.

The exceptional divisor above the singular point and the three geometrically
irreducible components of the intersection with the tangent cone constitute
an $I_4$ configuration, so we may form an elliptic fibration.  It is
defined by the equations $(x_2:x_3)$.  In addition to
the $I_4$ fibre, it also has an $I_2$ fibre consisting of two conics $C_2, C_3$
that meet in two points.  Let $C_2$ be the one that is not locally solvable
at $3, \infty$ and $C_3$ the one with rational points such as $(-2:-1:2:-1:1)$.

The curves $L_3,L_4,C_3$ now form an $I_3$ fibre.  The associated fibration
also has an $I_2$ fibre and another $I_3$, and is defined by
$(x_0+6x_2/5:x_1+x_2/3-x_3)$.  From the $I_3$ fibre that contains
$C_1$ we obtain a new conic $C_4$ which is the union of two lines defined
over $\Q(\sqrt{-3})$, while the $I_2$ fibre affords a curve $Q_1$ of degree $4$.

We now find some additional curves by mapping to another model in projective
space.  To be precise, consider the forms of degree $2$ that vanish on
$C_1, L_1, C_4$ mod those vanishing on $S_{29}$.  These give us a map to
a surface in $\P^3$, and pulling back the singular points and components of
their tangent cones we find three new conics $C_5, C_6, C_7$ of which
$C_5$ contains $(0:1:0:1:1)$, while $C_6$ is reducible over $\Q(\sqrt{-1})$
and $C_7$ contains $(-6/5:0:1:0:1)$.

The fibration with $C_2 \cup C_3 \cup C_5$ is the one we are looking for.
As in \cite[Section 14]{elkies-kumar}, it has a rational $I_4$ fibre, a
rational $I_3$ and three defined over the cubic field of discriminant
$-87$, and three $I_2$ fibres defined over the cubic field of discriminant
$-116$.  In addition, there is a rational fibre of type $II$.  It is easy
to change coordinates so that the $I_4, I_3, II$ fibres of the two are
above the same points of $\P^1$.  This only shows that the Jacobian of the
fibration we constructed on $S_{29}$ is isomorphic to the fibration of
\cite{elkies-kumar}, but $L_1$ is a section of the fibration, so the
two are isomorphic and we have found the birational equivalence.
The defining equations for this fibration (after changing coordinates) are
$$(x_0-3x_1/5-7x_3/5:x_2+2x_3).$$

\subsection{$\Q(\sqrt{37})$}
In this case it is a bit more difficult to get started, because we have no
singular points.  By searching we can easily find three lines on the
surface, defined by the following equations:
$$\begin{aligned}
&  L_1: x_0 = 9x_1 + x_4 = x_3 = x_5 = 0, \cr
&  L_2: x_0 = 9x_1 - x_4 = x_3 = x_5 = 0, \cr
&  L_3: 27x_0 + 2x_5 = 27x_1 - 3x_4 + 8x_5 = 3x_2 + x_4 - 2x_5 = 9x_3 + x_5 = 0, \cr
\end{aligned}$$
(and with only slightly more difficulty we could find three more).
Here $L_2$ meets $L_1, L_3$, but $L_1$ and $L_3$ are disjoint, so no elliptic
fibration has a fibre supported on $L_1 \cup L_2 \cup L_3$.
Letting $H$ be the hyperplane class, we easily
compute that $(H-[L_1]-[L_2]-[L_3])^2 = 0$, where $[L_i]$ is the Picard class
of $L_i$, and as this class is easily checked to be free of base components we
obtain an elliptic fibration from the linear forms $27x_0 + 2x_5, 9x_3 + x_5$
vanishing on all three lines.  The reducible fibres are four of type $I_2$,
two of which consist of curves of degree $1, 4$ and two $2,3$.  Let $C_1$
be the conic in a reducible fibre that contains the point $(-4:1:0:0:9:0)$.

We are then fortunate to find that $C_1, L_1, L_2$ all meet transversely at
$(0:0:1:0:0:0)$, thus forming a fibre of type $IV$.  As in
\cite[Section 16]{elkies-kumar}, the associated fibration turns out in addition
to have a rational $I_2$ fibre and six of type $I_3$, of which two are rational
and the other four are conjugate over a quadratic extension of $\Q(\sqrt{37})$.
There are two possible ways to match rational $I_3$ fibres on the two, only
one of which is correct.  As in the previous section, we find that the
Jacobians of the general fibres are isomorphic.  Since $L_3$ is a section,
this is indeed an elliptic fibration, and so in fact the two fibrations are
isomorphic as curves of genus $1$ over $\Q(t)$.  The defining polynomials are
$$(-x_0-2x_1-x_2+2_3:x_0+2x_1+x_2-x_3).$$

This completes the proof that the models of \cite{williams} are birationally
equivalent to those of \cite{elkies-kumar}.
\bibliography{art}

\end{document}